\documentclass[a4paper,reqno]{amsart}
\usepackage[utf8]{inputenc}
\usepackage{amsmath,amssymb,amsthm,amsfonts}
\usepackage{bbm}
\usepackage{euscript}
\usepackage{enumitem}
\usepackage{nicefrac}
\usepackage{mathtools}	

\newtheorem{theorem}{Theorem}
\newtheorem{lem}{Lemma}

\newtheorem{prop}{Proposition}

\theoremstyle{remark}
\newtheorem{remk}{Remark}
\newtheorem{exam}{Example}

\theoremstyle{definition}

\allowdisplaybreaks

\newcommand{\e}{\mathrm{e}}
\newcommand{\Leb}{\mathrm{Leb}}
\newcommand{\ov}{\overline}
\newcommand{\supp}{\mathop{\rm supp}}

\DeclareMathOperator\E{E}

\DeclareMathOperator\Prob{P}

\newcommand{\vk}{\varkappa}

\newcommand{\mbR}{{\mathbb R}}

\newcommand{\mbN}{{\mathbb N}}

\newcommand{\cD}{{\mathcal D}}
\newcommand{\cB}{{\mathcal B}}
\newcommand{\cE}{{\mathcal E}}

\newcommand{\wt}{\widetilde}

\newcommand{\pt}{\partial}
\newcommand{\1}{1\!\!\,{\rm I}}

\DeclareMathOperator*{\esssup}{ess\,sup}

\usepackage[	
	backend=biber, 
	style=alphabetic,
	bibencoding=utf8,
	autolang=other,
	maxcitenames=3,
	maxbibnames=10,
	giveninits=true,
	date=year,
	doi=false,isbn=false,url=false,
	eprint=true
	]{biblatex}

\DeclareFieldFormat[book,incollection]{series}{#1}
\renewbibmacro{in:}{}
\AtEveryBibitem{\clearlist{language}}
\AtEveryBibitem{\clearfield{note}}
\begin{document}

\title{On $1-$point densities for Arratia flows with drift}

\author[A.A.~Dorogovtsev]{Andrey A. Dorogovtsev}
\author[M.B.~Vovchanskyi]{Mykola B. Vovchanskyi}
\address{Institute of Mathematics,  National Academy of Sciences of Ukraine, Tereshchenkivska str. 3, 01601, Kiev, Ukraine}

\maketitle

\begin{abstract}
We show that if drift coefficients of Arratia flows converge in $L_1(\mathbb R)$ or $L_{\infty}(\mathbb R)$ then the 1-point densities associated with these flows converge  to the density for the flow with the limit drift.


\end{abstract}

The Arratia flow, a continual system of coalescing Wiener processes that are independent until they meet, was introduced independently as a limit of coalescing random walks in \cite{Ar79Coalescing}, as  a system of reflecting Wiener processes in \cite{SouBaWer00Reflection} and as a limit of stochastic homeomorphic flows in \cite{Do04OneBrownian}. If interpreted as a collection of particles started at $0,$ it is a part of the Brownian web \cite{FonIsoNewRa04Brownian}.  At the same time, one can construct the Arratia flow with drift  using flows of kernels defined in \cite{JanRai04Flows} (see \cite[{\S 6}]{Ria18Duality} for a short explanation) or directly via martingale problems by adapting the method used in \cite{Ha84Coalescing} to build coalescing stochastic flows with more general dependence between particles (see \cite[Chapter 7]{Do07MeasureEng} for this approach).
 Following \cite[Chapter 7]{Do07MeasureEng}, we consider a modification of the Arratia flow that introduces drift affecting the motion of a particle within the flow. So by an Arratia flow $X^a\equiv\{X^a(u,t)\mid u\in\mbR, t\in\mbR_+\}$ with bounded measurable drift $a$ we understand a collection of random variables such that 
\begin{enumerate}
\item for every $u$ the process $X^a(u,\cdot)$ is an It\^o process with diffusion coefficient $1$ and drift $a;$ 
\item for all $t\ge 0$ the mapping $X^a(\cdot,t)$ is monotonically increasing;
\item for any $u_1, u_2$ 
the joint quadratic covariation of the martingale parts of $X^a(u_1,\cdot)$ and $X^a(u_2,\cdot)$ equals $(t-\inf\{s\mid X^a(u_1,s)=X^a(u_2,s) \})_+,$ 
with $\inf \emptyset$ being equal $\infty$ by definition.  
\end{enumerate}

\begin{theorem}[{\cite[Lemma 7.3.1]{Do07MeasureEng}}]
For any measurable and locally bounded function $a$ the flow $X^a$ exists and is unique in the terms of finite-dimensional distributions.
\end{theorem}
\begin{remk}
\label{remark:existence}
\cite[Lemma 7.3.1]{Do07MeasureEng} is stated for Lipschitz continuous drift, but, being based on the classical existence and uniqueness theorems for the martingale problem, the proof holds verbatim for any measurable locally bounded drift.
\end{remk}   

Since the set $\{X^a(u,t)\mid u\in \mbR\}$ is known to be locally finite for all $t>0$ \cite[Chapter 7]{Do07MeasureEng}, one defines for any $t>0$ the point process $\{X^a(u, t) \mid u \in \mbR\}.$  
Studying of such a process can be performed using point densities (see: \cite{MunRajTriZab06Multi, TriZab11Pfaffian} for a definition and a representation in terms of Pfaffians in the case of zero drift, respectively; \cite{DorVov20ApproximationsEng, DoVov20Representations} for representations in the case of non-trivial drift; \cite{GliFom18Limit,Fo16Distribution} for applications to the study of the above-mentioned point process). In accordance with \cite[Appendix B]{MunRajTriZab06Multi}, 
the following definition is adopted in the present paper: the 1-point density at time $t$ of the point process associated with $X^a$ is given via
\begin{equation}
\label{eq:density.defin}
p^{a,1}_t(u) = \lim_{\delta\to 0+} \delta^{-1} \Prob \big( X^a(\mbR, t) \cap [u;u+\delta]  \not= \emptyset \big),
\end{equation}
where the limit exists a.e.

The well-known relation between PDEs and SDEs is used to show that $p^{a_n,1}_t\to p^{a,1}_t, n\to \infty, a.e.$  whenever $a_n\to a, n\to\infty,$ in $L_1(\mbR)$ or $L_\infty(\mbR)$ under additional natural assumptions, which is the main result of the article.
It should be noted that one may approach the problem stated utilizing the representation of $p_t^{a,1}$ given in \cite[Theorem 2.1]{DoVov20Representations}. However, the method used in the present paper allows us to avoid a discussion of properties (e.g. uniform integrability) of certain stochastic exponentials involving Brownian bridges. An additional benefit is that the term equal to $p^{0,1}_t$ is separated from those containing drift in the resulting expression. Thus there is a possibility to obtain local estimates of the density in the terms of drift,  in particular, in the studies of  point concentration when external force modeled with non-zero  drift is present. Moreover, following the route proposed in \cite[\S 4, Appendixes A, B]{MunRajTriZab06Multi} one can estimate $n-$dimensional densities in the terms of $p^{a,1}_t$ given estimates for the corresponding transitional density. 

\begin{prop}
\label{prop:density.existence}
The density $p^{a,1}_t$ exists. 
\end{prop}
\proof
Following \cite[Appendix B]{MunRajTriZab06Multi} one considers the measure $\nu^{a,1}$ on $\cB(\mbR)$ given via
\[
\nu_t^{a,1}(A) = \E |X^a(\mbR, t) \cap A|, \quad A \in \cB(\mbR),
\]
where given $B\in \cB(\mbR)$ $X^a(B, t) = \{X^a(u,t)\mid u\in B\}$  by definition and $|B|$ equals the number of distinct points in $B.$
Then 
\begin{equation}
\label{eq:measure.as.limit}
\E |X^a(\mbR, t) \cap A| = \lim_{U\to\infty} \E |X^a([-U;U], t) \cap A|.
\end{equation}

The Girsanov theorem for the Arratia flow \cite[Theorem 7.3.1]{Do07MeasureEng} states that given $T>0$  and a bounded Lipshitz continuous function $a$ the distribution of $X^a$ as a random element in the Skorokhod space $D(\mbR, C([0; T]))$ is absolutely continuous w.r.t the distribution of $X^0.$  To extend this result to the case of $a\in L_\infty(\mbR),$ one proceeds as follows. The result of \cite[{Lemma 7.3.2}]{Do07MeasureEng} states that for any $t>0$ there exists $C_t > 0$ such that
\begin{equation}
\label{eq:prob.not.meeting.estimate}
\Prob\left\{ X^a(u, t) \not= X^a(v, t) \right\} \le C_t (v-u), \quad  u, v\in\mbR, u<v.
\end{equation}
and is valid for $a\in L_\infty(\mbR).$ Using this estimate one  extends the proof of \cite[Theorem 2]{Do05Remarks}, so the flow $X^a$ admits a version in $D(\mbR, C([0; T])).$ Then it is left to notice that the proof of \cite[Theorem 7.3.1]{Do07MeasureEng} uses only the standard Girsanov theorem and the estimate \cite[{Lemma 7.3.2}]{Do07MeasureEng}. Therefore, the Girsanov theorem for the Arratia flow holds for $a\in L_\infty(\mbR),$ too.

Since $\nu_t^{0,1}$ is absolutely continuous w.r.t. the Lebesgue measure~\cite[Appendix B]{MunRajTriZab06Multi}, for any $A\in\cB(\mbR) $ with $\Leb(A)=0$ the assumption
\[
|X^0([-U;U], t) \cap A| = 0 \quad \mathrm{a.s.}, \quad U \ge 0,
\] 
implies
\[
\E |X^a([-U;U], t) \cap A | = 0, \quad U \ge 0.
\]
Hence by \eqref{eq:measure.as.limit} there exists the density $q_t^{a,1}$ such that 
\[
\int_A q_t^{a,1}(x) dx = \nu_t^{a,1}(A) = \E \sum_{v\in X^a(\mbR, t)} \1_A(v), \quad  A\in\cB(\mbR).
\]

Analogously one proves the existence of $n-$dimensional densities $q_t^{a,n}$ such that
\[
\E\sum_{
\begin{subarray}{c}v_1, \ldots, v_{n}\in X^a(\mbR, t),\\
v_1, \ldots, v_{n} \ \mbox{\small are distinct}
\end{subarray}
}
\1_A((v_1, \ldots, v_{n}))=\int_{A}q_t^{a,n}(x)dx = \nu_t^{a,n}(A), \quad A \in \cB(\mbR^n).
\]
By the Lebesgue differentiation theorem, for almost all $x=(x_1, \ldots, x_n)\in\mbR^n,$
\begin{equation*}
\label{eq:p.leb.diff}
q_t^{a,n}(x) = \lim_{\delta\to 0+} \delta^{-n} \nu_t^{a,n} \left( \mathop{\times}_{j=\ov{1,n}} [x_k;x_k+\delta] \right).
\end{equation*}
In particular,  for almost all $x$
\[
q_t^{a,1}(x) = \lim_{\delta\to 0+} \delta^{-1} \E |X^a(\mbR, t) \cap [x;x+\delta]|.
\]
In order to obtain \eqref{eq:density.defin}, one still needs to replace $|X^a(\mbR, t) \cap [x;x+\delta]|$ with $\1[X^a(\mbR, t) \cap [x;x+\delta] \not= \emptyset]$ in the last expression. The argument presented in  \cite[Appendix B]{MunRajTriZab06Multi} for $a=0$ is applicable here. 



If $a$ is additionally assumed to be Lipschitz continuous it is shown in \cite[Corollary 6.1]{Ria18Duality} that the dual flow \cite[Section 2.2]{TriZab11Pfaffian} to the flow $X^a$ is $X^{-a}.$ The same proof can be extended to the case of a just bounded function $a.$ Thus, due to the analog of \eqref{eq:prob.not.meeting.estimate} for  $-a$ we have, for $u,v\in\mbR, u < v,$
\begin{align*} 
\sup_{x\in\mbR} \Prob\left\{  X^a(x,t) \in [u;v]  \right\} &\le \Prob\left\{  X^a(\mbR,t) \cap [u;v] \not = \emptyset \right\}  \nonumber  \\
&= \Prob\left\{ X^{-a}(u,t) \not = X^{-a}(v,t) \right\} \nonumber \\
&\le C_t(v-u). 
\end{align*} 

Using this estimate and repeating the arguments in \cite[\S 4]{MunRajTriZab06Multi}, we get
\begin{align*}
\label{eq:induction.est.5}
&\Prob\left\{ X^a(\mbR, t) \cap [u_k,v_k] \not = \emptyset, k =\ov{1,n}\right\} \le C_t^n \prod_{k=\ov{1,n}} (v_k - u_k), \notag \\
& \qquad  u_1 < v_1 < u_2 < \ldots u_n < v_n, \notag  \\
& \qquad n \in \mbN. 
\end{align*} 
Exactly as in \cite[Appendix B]{MunRajTriZab06Multi} we then conclude
\begin{equation*}
\nu_t^{a,n} (A) \le C_t^n \, \Leb(A), \quad A\in \cB(\mbR^n), n\in\mbN. 
\end{equation*}
Thus  
\begin{align*}
\lim_{\delta \to 0+} & \delta^{-1} \E \left( |X^a(\mbR, t) \cap [x;x+\delta]| - \1\big[ X^a(\mbR, t) \cap [x;x+\delta] \not= \emptyset \big] \right) \\
& \quad \le \lim_{\delta \to 0+} \delta^{-1} \E |X^a(\mbR, t) \cap [x;x+\delta]| \big( |X^a(\mbR, t) \cap [x;x+\delta]| - 1\big) \\
& \quad = \lim_{\delta \to 0+} \delta^{-1} \nu_t^{a,2} \big([x;x+\delta] \times [x;x+\delta]\big) \\
& \quad \le  C_t^2 \lim_{\delta \to 0+} \delta \\ 
& \quad = 0,
\end{align*}
which concludes the proof. \qed

\begin{remk}
If $a$  is additionally Lipschitz continuous one can prove that $p_t^{a,n}$ has  a continuous version, using the aforementioned representations from \cite{DoVov20Representations} or the analog of \cite[Theorem 3.1]{DorVov20ApproximationsEng}. The same conclusion for $a\in L_\infty(\mbR)$ follows from Theorem \ref{th:density.prop.repr}.
\end{remk}

Let $D_2 = \{u\in\mbR^2\mid u_1 < u_2\},$ and put for $a\in L_\infty(\mbR)$
\[
\nabla^a_x = \sum_{k=1}^2 a(x_k) \pt_{x_k}, \quad x=(x_1,x_2) \in \mbR^2.
\]   
In what follows, we construct a series  $W^a \in C(\ov D_2 \times (0;\infty)) \cap C(D_2 \times [0;\infty)) \cap C^{1,0}(D_2 \times [0;\infty))$ such that  $W^a$ is a distribution solution \cite[p. 31]{Fol95Introduction} in $\cD^\prime(D_2\times (0;t)),$ the space of (Schwartz) distributions over infinitely differentiable compactly supported in $D_2\times (0;t)$ test functions, to 
\begin{align}
\label{eq:main.pde}
\pt_s W^a & = \frac{1}{2} \Delta W^a -  \nabla^a_x W^a.
\end{align} 
Additionally, as a function, $W^a$ satisfies
\begin{align}
\label{eq:main.pde.part2}
W^a(x,0) &= 1, \quad x \in D_2,\nonumber \\ 
W^a(x,s) &= 0,  \quad x\in \pt D_2, \ s > 0.
\end{align} 
The density $p^{a,1}_t$ admits a representation as a derivative of $W^{a}.$ As a result, in order to prove the convergence of the 1-point densities it is sufficient to consider he convergence of the corresponding series. 

It is supposed $a\in L_\infty(\mbR)$ throughout the paper. 
Let $x=(x_1,x_2)\in D_2$ and let $\xi^a_x=(\xi^a_{x_1},\xi^a_{x_2})$ be the unique weak solution of the Cauchy problem \cite[Corollary 3.11]{KaShre91Brownian} 
\begin{align}
\label{eq:xi^a}
d\xi^a_{x_k}(t)& = - a(\xi^a_{x_k}(t))dt + dw_k(t), \nonumber \\
\xi^a_{x_k}(0)& = x_k, \quad k=1, 2,
\end{align}
where $w_1, w_2$ are independent standard Wiener processes started at $0.$ Define
\begin{equation}
\label{eq:theta^a}
\theta^a_x = \inf\{s \mid \xi^a_x \in \pt D_2\}.
\end{equation}

\begin{prop}
\label{lem:density.via.prob}
Let $a\in L_\infty(\mbR), t>0.$ The density $p^{a,1}_t$ admits the representation
\[
p^{a,1}_t(u) = \lim_{\delta\to 0+} \delta^{-1} \Prob\left(\theta^{a}_{(u,u+\delta)} > t\right),
\] 
whenever the limit exits.
\end{prop}
\proof As in the proof of Proposition \ref{prop:density.existence}, one gets, using the notion of the dual flow, 
\begin{align*}
\Prob \big( X^a(\mbR, t) \cap [u;u+\delta]  \not= \emptyset \big) = & \Prob \big( X^{-a}(u+\delta, t) > X^{-a}(u, t)\big)  \\
= & \Prob\left(\theta^{a}_{(u,u+\delta)} > t\right)
\end{align*}
for all $\delta>0.$ 
\qed

\begin{exam}
Using the previous result one immediately gets  $p^{0,1}_t(u) =  \frac{1}{\sqrt{\pi t}}.$ Let $a_c(x) = cx, c\not= 0.$ The existence of the Arratia flow $X^{a_c}$ follows from Remark \ref{remark:existence}. Moreover, using the arguments from \cite[proof of Lemma 2.2]{DorVov20ApproximationsEng}, it is possible to establish the analog of Proposition \ref{lem:density.via.prob}.   Solving the SDE for $\xi^{a_c}_{u_2}-\xi^{a_c}_{u_1}$ explicitly and using the arguments from \cite[proof of Lemma 2.2]{DorVov20ApproximationsEng} one gets
\[
p^{a_c,1}_t(u) = \sqrt{\frac{2}{\pi}} \frac{|c|^{1/2}}{\psi(t,c)}, \quad u\in\mbR,
\]
where
\begin{align*}
\psi(t,c) =& \begin{cases}
  \left(\e^{2tc}-1\right)^{1/2}, & c >0,\ t >0,\\    
  \left(1-\e^{2tc}\right)^{1/2}, & c <0,\ t >0.   
\end{cases}
\end{align*} 
Thus $p^{0,1}_t = \lim_{c\to 0} p^{a_c,1}_t$ $a.e.$ in $\mbR.$
\end{exam}

Note that due to the domain $D_2$ being unbounded the problem \eqref{eq:main.pde}--\eqref{eq:main.pde.part2} does not admit, in general, a unique solution. However, as usually, one expects the unique bounded solution $W^a$ to have a probabilistic representation, that is, to satisfy $W^a(x, s)=\Prob\left(\theta^{a}_{x} > s\right).$ At the same time, in the case of  $a$ smooth enough the Duhamel principle \cite[\S 2.3.1.c]{Evans10Partial}  states 
\[
W^a(x,s) = \int_{D_2} dy_0 \ g_{s}(x,y_0) + \int_0^s dr \int_{D_2} dy \ g_{s-r}(x,y) \left( - \nabla^a_y W^a(y,r)\right), 
\] 
where 
\[
g_r(x, y) = \frac{1}{2\pi r}\left( \e^{-\frac{\|x-y\|^2}{2r}} - \e^{-\frac{\|x-y^\ast\|^2}{2r}}\right),
\]
with $y^\ast = (y_2,y_1),$ is 
the transition density of the 2-dimensional Wiener process killed when it reaches $\pt D_2.$ Using this formal relation as a starting point we define, for any $a\in L_\infty(\mbR),$ 
\begin{align*}
W^a(x,s)=&\sum_{n\ge 0} W^a_{n}(x,s),\\
W^a_{0}(x,s) =& \int_{D_2} dy_0 \ g_{s}(x,y_0), \\
W^a_{n}(x,s)=& (-1)^n \int_{\Delta_n(s)} dr_1 \ldots dr_n \int_{D_2^{n+1}} dy_0 \ldots dy_n \ g_{s-r_n}(x,y_n) \prod_{j={1}}^n \nabla^a_{y_j} g_{r_j-r_{j-1}}(y_j,y_{j-1}), \\
& n\ge 1,
\end{align*}
where hereinafter $\Delta_k(r) = \{u\in\mbR^k\mid 0 \le u_1 \le \ldots \le u_k \le r\}, r> 0, k \ge 1,$ so that
\[
W^a_{n}(x,s) = - \int_0^s dr_n \int_{D_2} dy_n \ g_{s-r_n}(x,y_n) \nabla^a_{y_n} W^a_{n-1}(y_n, r_n), \quad n \ge 1.
\]
Note that always 
\begin{align*}
W^a_{0}((x_1,x_2),s) &= \Prob\left(\theta^0_{(x_1,x_2)} > s\right) = \Prob\left( \sup_{r \le s} w(r) < \frac{x_2-x_1}{\sqrt{2}} \right) \\
&= \sqrt\frac{2}{\pi} \int_0^{\frac{x_2-x_1}{\sqrt{2s}}} da \ \e^{-\frac{a^2}{2}},
\end{align*}
where $w$ is a standard Wiener process, and
\begin{align*}
\pt_{x_2}W^a_0((u,u),s) = - \pt_{x_1}W^a_0((u,u),s) = \frac{1}{\sqrt{\pi s}}, \quad s > 0.
\end{align*}

The next lemma gathers technical results from calculus that are needed afterwards. In what follows $\ov D_2$ is the closure of $D_2,$ the derivatives $\pt_{x_k}, k=1,2$ on $\pt D_2$ are understood as one-sided derivatives and, for $y=(y_1,y_2),$
\[
G_a(y) = \min\Big\{ \sum_{j=1}^2 |a(y_{j})|, \|a\|_{L_\infty(\mbR)} \Big\}.
\]
\begin{lem}
\label{lem:calculus}
\phantom{a} \\
\begin{enumerate} 
\item $\exists K, \gamma > 0 \ \forall a\in L_\infty(\mbR)  \ \forall r > 0$
\begin{align*}
&\left| \pt_{y_{1,k}} g_r(y_1,y_2) \right| \le  K  r^{-3/2} \e^{-\gamma \frac{\|y_1-y_2\|^2}{r}} \ \ in \ \ov D_2, \quad k = 1,2, \ \\
&\left| \nabla^a_{y_1} g_r(y_1,y_2) \right| \le  K G_a(y_1) r^{-3/2} \e^{-\gamma \frac{\|y_1-y_2\|^2}{r}} \ \ a.e. \ in \ \ov D_2, 
\end{align*}
where $y_1 =(y_{1,1}, y_{1,2});$
\item $\forall \alpha > 0  \ \forall y\in \mbR^2$
\[
\int_{D_2} dy_0 \  \e^{-\alpha\|y_0 -y\|^2} \le \frac{\pi}{\alpha};
\]
\item put $r_0 = 0;$ then $\forall s > 0 \ \forall n \ge 1$
\begin{align*}
&\int_{\Delta_n(s)}dr_1 \ldots dr_n \ (s-r_n)^{-1/2} \prod_{j=1}^n (r_j-r_{j-1})^{-1/2}= \frac{\pi^{\frac{n+1}{2}} s^{\frac{n-1}{2}}}{\Gamma\left(\frac{n+1}{2}\right)}, \\
&\int_{\Delta_n(s)}dr_1 \ldots dr_n \ \prod_{j=1}^n (r_j-r_{j-1})^{-1/2}= \frac{2\pi^{\frac{n}{2}} s^{\frac{n}{2}}}{n\Gamma\left(\frac{n}{2}\right)};
\end{align*}
where here and hereinafter $\Gamma$ is the gamma-function;
\item 
for any $n\ge 1$ the derivatives $\pt_{x_k} W^a_n, k=1,2,$ exist on $\ov D_2\times (0;+\infty),$ and $W^a_n, \pt_{x_k} W^a_n \in C(\ov D_2\times (0;+\infty));$ additionally,
\begin{align*}
\sup_{x\in \ov D_2}\left| W^a_n(x,s)\right| &\le C_1 \frac{C_2^{n} \|a\|_{L_\infty(\mbR)}^{n} s^{\frac{n}{2}}}{\Gamma\left(\frac{n}{2}\right)},  \\ 
\sup_{x\in \ov D_2}\left| \pt_{x_k} W^a_n(x,s)\right| &\le C_1 \frac{C_2^{n} \|a\|_{L_\infty(\mbR)}^{n} s^{\frac{n-1}{2}}}{\Gamma\left(\frac{n}{2}\right)},  \quad s > 0, k=1,2, n\ge 1,
\end{align*}
for some absolute positive constants $C_1, C_2$ independent of $s, x, n$ and $a;$
\item for any $t>0$ the series for $W^a$ converges uniformly on $\ov D_2\times[0;t];$ the derivatives $\pt_{x_k} W^a, k=1,2,$ exist and are continuous on $\ov D_2\times (0;+\infty).$  
\end{enumerate}
\end{lem}
\proof 
Denote the coordinates of the point $y_2$ by $y_{2,1}, y_{2,2},$ and put $y_2^\ast = (y_{2,2}, y_{2,1}).$ Note that, for $k=1,2,$
\begin{align*}
|\pt_{y_{1,k}} g_r(y_1,y_2) | &\le \frac{1}{2\pi r^2} \left(  \e^{-\frac{\|y_1-y_2\|^2}{2r}} |y_{1,k} - y_{2,k}| + \e^{-\frac{\|y_1-y_2^\ast\|^2}{2r}}  |y_{1,k} - y_{2,k}^\ast| \right)\\
& \le \frac{1}{2\pi r^2} \left(  \e^{-\frac{\|y_1-y_2\|^2}{2r}}  \|y_{1} - y_{2}\| + \e^{-\frac{\|y_1-y_2^\ast\|^2}{2r}}  \|y_{1} - y_{2}^\ast\| \right) \\
& \le \frac{C}{r^{3/2}} \left(   \e^{-\frac{\|y_1-y_2\|^2}{4r}} +  \e^{-\frac{\|y_1-y_2^\ast\|^2}{4r}} \right),
\end{align*}
for some absolute constant $C.$
Since $\e^{-\frac{\|y_1-y_2^\ast\|^2}{4r}} \le \e^{-\frac{\|y_1-y_2\|^2}{4r}}$ for $y_1,y_2 \in \ov D_2$ this yields Item (1).

 For Item (2), write 
\begin{align*}
\int_{D_2} dy_0 \  \e^{-\alpha\|y_0 -y\|^2} \le 2\pi \int_0^\infty d\rho \ \rho \e^{-\alpha\rho^2} = \frac{\pi}{\alpha}. 
\end{align*} 

Item (3) is a direct consequence of the integral representation given in \cite[\S 4.3]{Carl77Special} after the change of variables $r_1 = s u_1 , r_2 = s(u_1 + u_2), \ldots, r_n = s(u_1 + \ldots + u_n )$.\footnote{The authors are grateful to a colleague for providing the reference.} 


Fix $n$ and put $r_{n+1}=s, y_{n+1}=x.$ Combining Items (1) and (2) we get
\begin{align*}
\left| W^a_n(x,s)\right| &\le K^n \|a\|_{L_\infty(\mbR)}^{n}  \int_{\Delta_n(s)} dr_1 \ldots dr_n \int_{D_2^{n+1}} dy_0 \ldots dy_n  \ g_{s-r_n}(x,y_n) \times \\
& \qquad \qquad \qquad  \times \prod_{j=1}^{n}  \frac{\e^{-\gamma \frac{\|y_j-y_{j-1}\|^2}{r_j-r_{j-1}}}}{(r_j-r_{j-1})^{3/2}} \\
&\le K^n \|a\|_{L_\infty(\mbR)}^{n}  \left(\frac{\pi}{\gamma}\right)^{n} \int_{\Delta_n(s)} dr_1 \ldots dr_n \ \prod_{j=1}^{n} (r_j-r_{j-1})^{-1/2}, \\
\left|\pt_{x_k} W^a_n(x,s)\right| &\le K^{n+1} \|a\|_{L_\infty(\mbR)}^{n}  \int_{\Delta_n(s)} dr_1 \ldots dr_n \int_{D_2^{n+1}} dy_0 \ldots dy_n  \ \prod_{j=1}^{n+1}  \frac{\e^{-\gamma \frac{\|y_j-y_{j-1}\|^2}{r_j-r_{j-1}}}}{(r_j-r_{j-1})^{3/2}} \\
&\le K^{n+1} \|a\|_{L_\infty(\mbR)}^{n} \left(\frac{\pi}{\gamma}\right)^{n+1} \int_{\Delta_n(s)} dr_1 \ldots dr_n \ \prod_{j=1}^{n+1} (r_j-r_{j-1})^{-1/2}, 
\end{align*}
so the rest follows by calculus.
\qed

In \cite{PaWin08First} a general version of the problem \eqref{eq:main.pde}--\eqref{eq:main.pde.part2} in a bounded domain is studied in terms of weak (variational) solutions\footnote{The statement of \cite[Theorem 2.7]{PaWin08First} contains a misprint: $Q(t,x)$ should be equal $\Prob(\tau_A >  t).$} so the smoothness of the solution follows by standard results. Since it is easy to see that neither $\pt_{x_k x_k} W^a_n$ nor $\pt_{s} W^a_n$ exists for $n\ge 1$ this is not a viable approach for us, and we use the hypoellipticity of an operator $\frac{1}{2} \Delta -  \nabla^a_x$ for $a\in C^{\infty}(\mbR)\cap L_\infty(\mbR)$ instead to check that the  distribution solution $W^a$ is sufficiently smooth.

The proof of the following lemma is trivial. 
\begin{lem}
\label{lem:approx.a}
Let $M\in \cB({\mbR^m})$ for some $m\in\mbN,$ and $k\in L_1(M).$ Suppose functions  $f_n\colon M \mapsto \mbR, n\ge 0,$ are such that 
\[
\sup_{n\ge 0} \|f_n\|_{L_\infty(M)} < \infty.
\]
Suppose that one of the following conditions holds:
\begin{enumerate}
\item $f_n\in L_1(M), n\ge 0,$  and $f_n\to f_0, n\to \infty,$ in $L_1(M);$ 
\item  
$f_n\to f_0, n\to \infty,$ in $L_\infty(M).$ 
\end{enumerate}
Then 
\[
\int_M dy\ f_n(y) k(y)\to \int_M dy\ f_0(y) k(y), \quad n\to\infty.
\]
\end{lem} 
\begin{prop}
\label{prop:distr.solution.W_n}
For all  $n\ge 1$ in the sense of Schwartz distributions
\[
\pt_s W^a_n = \frac{1}{2} \Delta W^a_n-  \nabla^a_x W^a_{n-1}
\]
in $D_2 \times (0;\infty) .$
\end{prop}
\proof Suppose that $t >0$ is fixed throughout the proof. We have, for $n\ge 1,$
\[
W^a_n(x,s) = \int_0^s \int_{D_2} dr dy\ g_{s-r}(x, y) f_n(r, y), \quad x \in D_2, s \in (0;t),
\]
where
\[
f_n(r, y) = - \nabla^a_{y} W^a_{n-1}(y, r).
\]
By Item 4 of Lemma  \ref{lem:calculus}
\begin{equation}
\label{eq:bound.f_n}
\sup_{n\ge 1} \esssup_{r\in(0;t), y\in D_2} |f_n(r,y)| \le F r^{-\frac{1}{2}},
\end{equation}
for an absolute constant $F.$
Consider for $h>0$
\begin{align}
\label{eq:defn.H_k}
W^a_n(x, s+h) - W^a_n(x, s)&= \int_0^s \int_{D_2} dr dy\ \left( g_{s+h-r}(x,y) - g_{s-r}(x,y)\right) f_n(r,y)  \notag \\
& \quad + \int_s^{s+h} \int_{D_2} dr dy\  g_{s+h-r}(x,y) f_n(r,y) \notag \\
=& H_1(h,x,s) + H_2(h,x,s).
\end{align}
Given a compactly supported in $D_2$ function $v\in C^{\infty}(D_2)$  we have, by the Fubini theorem, properties of the Gaussian density and the Green formula,
\begin{align*}
\label{eq:10}
h^{-1}\int_{D_2} dx\ & v(x) H_1(h,x, s) = \int_{D_2} dx\ v(x)  h^{-1} \int_0^s \int_{D_2} dr dy \int_{s-r}^{s+h-r} d\tau\ \pt_\tau g_\tau(x,y)f_n(r,y) \notag \\
&= \int_0^s \int_{D_2} dr dy \int_{s-r}^{s+h-r} d\tau\ h^{-1} f_n(r,y) \int_{D_2} dx\ \pt_\tau g_\tau(x,y) v(x) \notag \\
&=  \int_0^s \int_{D_2} dr dy \int_{s-r}^{s+h-r} d\tau\ h^{-1} f_n(r,y) \int_{D_2} dx\ \frac{1}{2}\Delta_{x} g_\tau(x,y) v(x) \notag \\
&= -\frac{1}{2} \int_0^s \int_{D_2} dr dy \int_{s-r}^{s+h-r} d\tau\ h^{-1} f_n(r,y) \int_{D_2} dx\ \nabla_{x} g_\tau(x,y) \cdot \nabla_x v(x) \notag \\
&= -\frac{1}{2} \int_{D_2} \int_0^s \int_{D_2} dx dr dy \ \nabla_x v(x) \cdot \left( h^{-1} \int_{r}^{r+h} d\tau\  f_n(s-r,y)  \nabla_{x} g_\tau(x,y) \right), 
\end{align*}
where $\cdot$ stands for the inner product in $\mbR^2.$  For all $r$
\begin{equation*}
\label{eq:11}
h^{-1} \int_{r}^{r+h} d\tau\ \pt_{x_k} g_\tau(x,y) \to \pt_{x_k} g_{r}(x,y),  \quad k =1,2. 
\end{equation*}
Thus to prove  
\begin{align}
\label{eq:13}
h^{-1}\int_{D_2} dx\  v(x) H_1(h,x, s) &\to -\frac{1}{2}\int_{D_2}\int_0^s \int_{D_2} dx dr dy\  \nabla_x v(x) \cdot  f_n(r,y) \nabla_{x} g_{s-r}(x,y) \notag \\
&= -\frac{1}{2}\int_{D_2}dx \ \nabla_x v(x) \cdot \nabla_x W^a_n(x,s), \quad h\to 0+,
\end{align}
by the means of dominated convergence theorem it is sufficient to show that  both families
\begin{align*}
&\left\{  (x,r,y) \mapsto h^{-1} \pt_{x_k} v(x) f_n(s-r,y) \int_{r}^{r+h} d\tau\ \pt_{x_k} g_\tau(x,y) \mid h\in (0;1) \right\}, \\ &\qquad \qquad k= 1,2,
\end{align*}
are uniformly integrable on $\supp v \times (0;t) \times D_2.$ For that, note that by Lemma \ref{lem:calculus} and \eqref{eq:bound.f_n} for almost all $(r,y)$
\begin{align*}
\left| h^{-1} f_n(s-r,y) \int_{r}^{r+h} d\tau\ \pt_{x_k} g_\tau(x,y) \right| &\le KF h^{-1} (s-r)^{-\frac{1}{2}}\int_{r}^{r+h} d\tau\ \tau^{-\frac{3}{2}} \e^{-\gamma\frac{\|x-y\|^2}{\tau}}, 
\end{align*}
where, for some positive $C,$ 
\begin{align*}
h^{-1} \int_{r}^{r+h} d\tau\ \tau^{-\frac{3}{2}} \e^{-\gamma\frac{\|x-y\|^2}{\tau}} & \le 
\begin{cases}
   r^{-\frac{3}{2}} \e^{-\gamma\frac{\|x-y\|^2}{r}}, & \frac{2}{3} < \frac{r}{\gamma \|x-y\|^2}, \\
  (r+h)^{-\frac{3}{2}} \e^{-\gamma\frac{\|x-y\|^2}{r+h}}, & \frac{2}{3} > \frac{r+h}{\gamma \|x-y\|^2}, \\  
  \frac{C}{\|x-y\|^3}, & \frac{2}{3}\in \left[\frac{r}{\gamma \|x-y\|^2};\frac{r+h}{\gamma \|x-y\|^2}\right], \\    
\end{cases} 
\end{align*}
so that, for some $C_1, C_2 > 0,$   all  $\alpha\in(1;\frac{4}{3})$ and $\delta\in(0;1),$ we have, due to Item 2 of Lemma \ref{lem:calculus}, for almost all $r$
\begin{align*}
\label{eq:un.integr.est}
&  \int_{D_2} dy \left| f_n(s-r,y) h^{-1} \int_{r}^{r+h} d\tau\  \pt_{x_k} g_\tau(x,y) \right|^\alpha   \notag \\
& \qquad \le \frac{C_1}{(s-r)^{\frac{\alpha}{2}}} \int_{D_2} dy \Bigg\{ r^{-\frac{3\alpha}{2}} \e^{-\frac{\alpha\gamma \|x-y\|^2}{r}} + (r+h)^{-\frac{3\alpha}{2}} \e^{-\frac{\alpha\gamma \|x-y\|^2}{r+h}} \notag \\
& \qquad \qquad + \frac{1}{\|x-y\|^{3\alpha}} \1\left[\frac{2}{3}\in \left[\frac{r}{\gamma \|x-y\|^2};\frac{r+h}{\gamma \|x-y\|^2}\right]\right] \Bigg\} \notag \\
& \qquad \le \frac{C_1}{(s-r)^{\frac{\alpha}{2}}}  \Bigg\{  \frac{\pi}{\alpha\gamma}  \left( r^{1-\frac{3\alpha}{2}} + (r+h)^{1-\frac{3\alpha}{2}} \right) + 2\pi  \int_0^\infty d\rho \ \rho^{1-3\alpha}  \1\left[ \frac{r}{\rho^2} \le \frac{2\gamma}{3} \le \frac{r+h}{\rho^2} \right] \Bigg\} \notag \\
& \qquad \le \frac{C_1 \left(\frac{\pi}{\alpha\gamma}  + 2\pi\right)}{(s-r)^{\frac{\alpha}{2}}} \Bigg\{  \left( r^{1-\frac{3\alpha}{2}} + (r+h)^{1-\frac{3\alpha}{2}} \right) + \int_{\sqrt{ \frac{3r}{2\gamma}}}^{\sqrt{\frac{3(s+1)}{2\gamma}}} d\rho \ \rho^{1-3\alpha} \Bigg\} \notag \\
& \qquad \le  \frac{C_1 \left(\frac{\pi}{\alpha\gamma}  + 2\pi\right)}{(s-r)^{\frac{\alpha}{2}}}  \Bigg\{ \left( r^{1-\frac{3\alpha}{2}} + (r+h)^{1-\frac{3\alpha}{2}} \right)  + \left(\frac{3r}{2\gamma}\right)^{1-\frac{3\alpha+\delta}{2}} \int_0^{\sqrt{\frac{3(s+1)}{2\gamma}}} d\rho \ \rho^{-1+\delta}  \Bigg\} \notag \\
& \qquad \le \frac{C_2}{(s-r)^{\frac{\alpha}{2}}}  \Bigg\{ \left( r^{1-\frac{3\alpha}{2}} + (r+h)^{1-\frac{3\alpha}{2}} \right)  + r^{1-\frac{3\alpha+\delta}{2}} (s+1)^{\frac{\delta}{2}}  \Bigg\}.
\end{align*}
Choosing $\delta$ in such a way that  $1-\frac{3\alpha+\delta}{2} > -1$ gives, for some $C_3 >0,$ 
\begin{align*}
&\sup_{x\in D_2}  \int_0^s \int_{D_2} dr dy \left| h^{-1} f_n(s-r,y)  \int_{r}^{r+h} d\tau\ \pt_{x_k} g_\tau(x,y) \right|^\alpha   \notag \\
& \qquad \le C_2 \Bigg\{ 2 \int_{0}^{s}dr \ \frac{r^{1-\frac{3\alpha}{2}}}{(s-r)^{\frac{\alpha}{2}}} + (s+1)^{\frac{\delta}{2}} \int_0^s dr \ \frac{r^{1-\frac{3\alpha+\delta}{2}}}{(s-r)^{\frac{\alpha}{2}}}  \Bigg\} \notag \\
& \qquad \le C_3 \int_0^s dr \ \frac{r^{1-\frac{3\alpha+\delta}{2}}}{(s-r)^{\frac{\alpha}{2}}}  \notag \\
& \qquad < +\infty,
\end{align*} 
 
which yields \eqref{eq:13}.

For $H_2,$ we have, in terms of the process $\xi^a_x$  and the stopping time $\theta^a_x$ defined in \eqref{eq:xi^a}-\eqref{eq:theta^a}
\begin{align*}
h^{-1} H_2(h,x, s) = h^{-1} \int_s^{s+h} dr\ k_{x,s}(h,r),
\end{align*}
where
\[
k_{x,s}(h,r) = \E f_n\left(r, \xi^a_x(s+h-r)\right) \1\left[ \theta^a_x > s+h-r \right].
\]
Since $a$ is bounded, the Novikov condition holds on $(0;t)$ so by the Girsanov theorem 
\label{eq:17}
\begin{align*}
k_{x,s}(h,r) = \E f_n\left(r, \xi^0_x(s+h-r)\right) \1\left[ \theta^0_x > s+h-r \right] \cE_t^a,
\end{align*}  
where $\cE_t^a$ is the corresponding stochastic exponential. The family of random variables
\[
\left\{ f_n\left(r, \xi^0_x(s+h-r)\right) \1\left[ \theta^0_x > s+h-r \right] \cE_t^a \mid (r,h) \in A\right\},
\]
where $A = \{(u_1,u_2)\mid u_1\in [s;s+u_2], u_2\in [0;1]\},$ is uniformly integrable, so since the distribution of $\theta_x^0$ has no atoms the function $k_{x,s}$ is uniformly continuous on $A$ for any $(x,s).$ Fix $\epsilon > 0.$ Then there exists $h_0$ such that 
\[
h^{-1} \int_s^{s+h} dr\ |k_{x,s}(h_0,r) -k_{x,s}(h,r)| \le \epsilon
\]
whenever $h\le h_0.$ By the Lebesgue differentiation theorem for fixed $h_0$
\[
h^{-1} \int_s^{s+h} dr\ k_{x,s}(h_0,r) \to k_{x, s}(h_0,s), \quad h\to0+,
\]
where 
\[
k_{x, s}(h_0,s) \to k_{x, s}(0,s) = f_n(s,x), \quad h_0\to 0+.
\]
Therefore 
\begin{equation*}
\lim_{h\to 0+} h^{-1} H_2(h,x, s) = \lim_{h_0\to 0+} k_{x,s}(h_0,s) = f_n(s,x) = -\nabla^a_x W^a_{n-1}(x,s).
\end{equation*}

Consider the analog of \eqref{eq:defn.H_k} for a negative increment: for $h>0$
\begin{align*}
W^a_n(x, s-h) - W^a_n(x, s)&= \int_0^{s-h} \int_{D_2} dr dy\ \left( g_{s-h-r}(x,y) - g_{s-r}(x,y)\right) f_n(r,y)  \\
& \quad + \int_{s-h}^{s} \int_{D_2} dr dy\  g_{s-h-r}(x,y) f_n(r,y) \notag \\
=& \wt H_1(h,x,s) + \wt H_2(h,x,s),
\end{align*}
and the reasoning for $\wt H_1$ and $\wt H_2$ is the same as for their counterparts $H_1$ and $H_2,$ which concludes the proof.
\qed

\begin{theorem}
\label{th:W^a.prob.repr}
Let $a\in L_\infty(\mbR).$ For all $s>0$ and $x\in D_2$
\[
W^a(x,s) = \Prob\left( \theta^a_x > s\right)
\]
where $\theta^a_x$ is defined in \eqref{eq:theta^a}.
\end{theorem}
\proof Assume $a\in C^\infty(\mbR)$ additionally. It follows from Proposition \ref{prop:distr.solution.W_n} and the completeness of the space of Schwartz distributions  that in the sense of Schwartz distributions
\[
\pt_s W^a = \frac{1}{2} \Delta W^a-  \nabla^a_x W^a
\]
in $(0;\infty)\times D_2.$ Since the operator $\frac{1}{2} \Delta -  \nabla^a_x$ is hypoelliptic, $W^a\in C^\infty(D_2 \times (0;\infty))$ \cite[Theorem 3.4.1]{Stroock2008Partial}. By Item 4 of Lemma \ref{lem:calculus}, $W^a\in C(\ov D_2 \times (0;\infty)) \cap C(D_2 \times [0;\infty)).$ Since 
\[\
g_r((u,u), y) = 0, \quad r\ge 0, u\in\mbR, y\in\mbR^2,
\]   $W^a$ satisfies the Dirichlet boundary condition as a function. To show 
\begin{equation}
\label{eq:prob.repr}
W^a(x,s) = \Prob\left( \theta^a_x > s\right),
\end{equation}
one proceeds as follows.

Let
\[
K_n = \left\{ x\in D_2 \mid \|x\| \le n,  \inf_{y\in\pt D_2} \|x -y\| \ge \frac{1}{n} \right\}, \quad n\ge 1, 
\] 
and define, for fixed $T>0,$ 
\begin{align*}
\theta^{n}_{t,x; T} &= T \wedge \inf\left\{s \ge t \mid \xi^a_{t,x}(s) \not\in K_n \right\},  \\
\theta_{t,x} &=  \inf\left\{s \ge t \mid \xi^a_{t,x}(s) \not\in D_2 \right\}, \quad t > 0, x \in D_2,
\end{align*}
where, analogously to \eqref{eq:xi^a}, for a standard Wiener process $w$ 
\begin{align*}
d\xi^a_{t,x}(s)& = - a(\xi^a_{t,x}(s))ds + dw(s), \quad s \ge t, \nonumber \\
\xi^a_{t,x}(t)& = x.
\end{align*}
The standard reasoning via the It\'o formula and the exhaustion method (e.g. \cite[\S 2.2]{Freid85Functional} or \cite[Ch.8, \S 5]{GikhSko77IntroEng}) gives, for any fixed $h>0, t \in (0; T),$
\begin{equation}
\label{eq:repr.}
\E W^a\left(h+T-\theta^n_{t,x; T}, \xi^a_{t,x}(\theta^{n}_{t,x; T})\right) = W^a(h+T-t,x), \quad n\ge 1, x\in D_2.
\end{equation}
Since $W^a\in C(D_2 \times [0;\infty))$
\begin{equation*}
\label{eq:repr.rhs}
W^a(h+T-t,x) \to W^a(T-t, x), \quad h\to 0+.
\end{equation*}

We have 
\begin{align*}
\label{eq:repr.lhs}
\E W^a\left(h+T-\theta^n_{t,x; T}, \xi^a_{t,x}(\theta^{n}_{t,x; T})\right) = \E W^a\left(h + T - t - \theta^n_{0,x; T-t}, \xi^a_{0,x}(\theta^n_{0,x; T-t})\right).
\end{align*}
The process $\xi^a_{0,x}$ is bounded on $[0;T-t]$ with probability $1.$ On the set $\{\theta_{0,x} \le T-t\}$
\[
W^a\left(h + T - t - \theta^n_{0,x; T-t}, \xi^a_{0,x}(\theta^n_{0,x; T-t})\right)  \to W^a\left(h + T - t - \theta_{0,x}, \xi^a_{0,x}(\theta_{0,x})\right) = 0, \quad n\to\infty,
\]
since $W^a\in C(\ov D_2 \times (0;\infty)).$ On the set $\{\theta_{0,x} > T-t\}$
\[
W^a\left(h + T - t - \theta^n_{0,x; T-t}, \xi^a_{0,x}(\theta^n_{0,x; T-t})\right) \to W^a\left(h, \xi^a_{0,x}(T-t)\right), \quad n\to\infty,
\]
and, since $\xi^a_{0,x}(T-t)\in D_2$ now and $W^a\in C(D_2 \times [0;\infty)),$
\[
W^a\left(h, \xi^a_{0,x}(T-t)\right) \to W^a(0,x)=1, \quad h\to 0+.
\]
Returning to \eqref{eq:repr.}, passing to limit w.r.t $n$ and $h$ and introducing $s= T-t > 0$, we get
\[
W^a(x, s) =  \E \1\left[\theta_{0,x} > s \right], 
\]
which coincides with \eqref{eq:prob.repr}.

Consider the case $a\in L_{\infty}.$ Let $\eta$ be a mollifier and put $a_n = a \ast n\eta(\frac{\cdot}{n}), n\ge 1.$ Additionally put $a_0 = a.$ Then $a_n\in  C^\infty(\mbR),$ $|a_n| \le |a|, n\ge 1,$ and $a_n\to a, n\to\infty,$ in $L_1(\mbR).$ Consider $s\in (0;t)$ for a fixed $t.$ Then the Novikov condition
\[
\sup_{n\ge 0} \E \exp\Big\{ \frac{1}{2} \int_0^t ds \sum_{j=1}^2 a_n^2 \left(\xi^a_{x_j}(s)\right)\Big\} < \infty,
\]
holds, where $\xi^a_x = \xi^a_{0, x}.$ Thus, by the Girsanov theorem, 
\begin{align*}
\Prob\left(\theta^{a_n}_x > s\right) =& \E \vk\cE_n,
\end{align*}
where 
\begin{align*}
\vk =& \1\left[\theta^0_x > s\right], \\
\cE_n =& \exp\Big\{\sum_{j=1}^2 \int_0^t  a_n \left(\xi^0_{x_j}(r)\right) d\xi^0_{x_j}(r) - \frac{1}{2}\sum_{j=1}^2 \int_0^t dr\ a_n^2 \left(\xi^0_{x_j}(r)\right)  \Big\}, \\
& n\ge 0.
\end{align*}
Since $a\in L_\infty(\mbR)$ the Cauchy inequality implies that the sequence $\{\vk\cE_n\}_{n\ge 1}$ is uniformly integrable. Consequently, in order to prove 
\begin{equation}
\label{eq:add.1}
\Prob\left(\theta^{a_n}_x > s\right) = \E \vk\cE_n \to  \E\vk\cE_0 = \Prob\left(\theta^{a}_x > s\right), \quad n\to\infty,
\end{equation} 
it is sufficient to check that $\cE_n\to \cE_0$ in probability. For that, consider, for fixed $j,$
\begin{align}
\label{eq:add}
\E &\left|\int_0^t  dr\ a_n^2 \left(\xi^0_{x_j}(r)\right) - \int_0^t  dr\ a_0^2 \left(\xi^0_{x_j}(r)\right) \right|\le \int_0^t dr \E \left| a_n^2 \left(\xi^0_{x_j}(r)\right) -a_0^2 \left(\xi^0_{x_j}(r)\right) \right|.
\end{align}
Since for every $s$  the sequence
\[
\Big\{ u\mapsto  \frac{\left| a^2_n(u) -a^2_0(u)  \right|\e^{-\frac{\|u-x_j\|^2}{2s}}}{s} \mid n\ge 1 \Big\}
\] 
is uniformly integrable w.r.t. the Lebesgue measure in $\mbR$  the dominated convergence theorem implies
\[
\E \left| a_n^2 \left(\xi^0_{x_j}(r)\right) -a_0^2 \left(\xi^0_{x_j}(r)\right)\right| \to 0, \quad n\to\infty,
\]
for all $s\in(0;t).$ Using the dominated convergence theorem one more time in \eqref{eq:add} we obtain
\[
\int_0^t  dr\ a_n^2 \left(\xi^0_{x_j}(r)\right) \to \int_0^t  dr\ a_0^2 \left(\xi^0_{x_j}(r)\right), \quad n\to\infty,
\]
in probability. Analogously, one proves 
\[
\int_0^t  a_n \left(\xi^0_{x_j}(r)\right) d\xi^0_{x_j}(r) \to \int_0^t  a_0 \left(\xi^0_{x_j}(r)\right) d\xi^0_{x_j}(r), \quad n\to\infty,
\]
in probability, which finishes the proof of \eqref{eq:add.1}.

Thus
\[
W^{a_n}(x,s) = \Prob\left( \theta^{a_n}_x > s\right) \to \Prob\left( \theta^{a}_x > s\right), \quad n\to\infty.
\]
Therefore, it is left to show 
\begin{equation}
\label{eq:add.4}
W^{a_n}(x,s)\to W^{a}(x,s), \quad n\to\infty.
\end{equation}
Note that by Item 4 of Lemma \ref{lem:calculus} for all $s\in(0;t)$ and $x\in D_2$
\begin{equation}
\label{eq:add.2}
\sup_{n\ge 1}  \left| \sum_{m\ge M} W^{a_n}_m(x,s) -\sum_{m\ge M} W^{a}_m(x,s)\right| \to 0, \quad M \to\infty.
\end{equation}
For all $m\ge 1$ by Lemma \ref{lem:calculus}
\begin{align*}
&\Big| \prod_{j={1}}^m \nabla^a_{y_j} g_{r_j-r_{j-1}}(y_j,y_{j-1}) - \prod_{j={1}}^m \nabla^{a_n}_{y_j} g_{r_j-r_{j-1}}(y_j,y_{j-1})\Big| \\
& \qquad \le \sum_{k=1}^m \Big| \prod_{j={1}}^{k-1} \nabla^{a_n}_{y_j} g_{r_j-r_{j-1}}(y_j,y_{j-1})\  \nabla^{a_n-a}_{y_k} g_{r_k-r_{k-1}}(y_k,y_{k-1}) \prod_{j={k+1}}^m \nabla^{a}_{y_j} g_{r_j-r_{j-1}}(y_j,y_{j-1})\Big| \\
&\qquad \le K^m \|a\|_{L_\infty(\mbR)}^{m-1} \prod_{j={1}}^{m}  \frac{\e^{-\gamma \frac{\|y_1-y_2\|^2}{r_j-r_{j-1}}}}{(r_j-r_{j-1})^{3/2}}  \sum_{k=1}^m  \sum_{j=1}^2 |a_n(y_{k,j})-a(y_{k,j})| ,
\end{align*}
where the product over the empty set is equal to $1$ by definition. Therefore, due to \eqref{eq:add.2} and the definition of $W^a$ it is sufficient to show  that
\[
\int_{\Delta_m(s)} dr_1 \ldots dr_m \int_{D_2^{m+1}} dy_0 \ldots dy_m\ \prod_{j={1}}^{m}  \frac{\e^{-\gamma \frac{\|y_1-y_2\|^2}{r_j-r_{j-1}}}}{(r_j-r_{j-1})^{3/2}} |a_n(y_{k,j})-a(y_{k,j})| \to 0, \quad n\to\infty,
\]
for fixed $m$ and $j.$ Setting $f_n(y_{k,j}) = |a_n(y_{k,j})-a(y_{k,j})|$ in Lemma \ref{lem:approx.a} and applying Items 2 and 3 of Lemma \ref{lem:calculus} finishes the proof. We omit the details.\qed

\begin{theorem}
\label{th:density.prop.repr}
Let $a\in L_\infty(\mbR).$ Then for any $t>0$ the density $p^{a,1}_t$ has a continuous version
\[
p^{a,1}_t(x) = \pt_{x_2} W^a((x,x),t) = \sum_{n\ge 0} \pt_{x_2} W^a_n((x,x),t).
\]
\end{theorem}
\proof Combining Proposition \ref{lem:density.via.prob},  Theorem \ref{th:W^a.prob.repr} and the fact that $W^a((x,x), t) = 0$ we get 
\begin{align*}
p^{a,1}_t(x) &= \lim_{\delta\to 0+} \delta^{-1} \Prob\left(\theta^{a}_{(x,x+\delta)} > t\right) \notag \\
&= \lim_{\delta\to 0+} \delta^{-1} \left(W^a((x,x+\delta),t) - W^a((x,x), t)\right),
\end{align*}
so the conclusion follows by Lemma \ref{lem:calculus} and calculus.
\qed

The next result is stated for a continuous version of the $1-$dimensional density. 
\begin{theorem}
\label{th:density.conver}
Assume $a_n\in L_\infty(\mbR), n\ge 0; \sup_{n\ge 0} \|a_n\|_{L_\infty(\mbR)}< \infty.$ Let one of the following conditions hold:
\begin{enumerate}
\item $a_n\to a_0, n\to\infty,$ in $L_\infty(\mbR);$ 
\item $a_n\in L_1(\mbR), n\ge 0$, and $a_n\to a_0, n\to\infty,$ in $L_1(\mbR).$ 
\end{enumerate}
Then for any $t>0$ and all $x\in\mbR$
\[
p^{a_n,1}_t(x)\to p^{a_0,1}_t(x), \quad n\to\infty.
\]
\end{theorem}
\proof By Theorem \ref{th:density.prop.repr}, the conclusion of the theorem is equivalent to 
\[
\sum_{k\ge 0} \pt_{x_2} W^{a_n}_k((x,x),t) \to \sum_{k\ge 0} \pt_{x_2} W^{a_0}_k((x,x),t), \quad n\to\infty.
\]
Since $\sup_{n\ge 0} \|a_n\|_{L_\infty(\mbR)}< \infty,$ one applies the same reasoning as in the proof of \eqref{eq:add.4}, utilizing the uniform  estimates for $\pt_{x_k}W^{a_n}(x,s)$ of Item 4 in Lemma \ref{lem:calculus}, and uses both statements of Lemma \ref{lem:approx.a} subsequently. We omit the details. \qed

\section*{Acknowledgements}
The authors are grateful to the anonymous referee for the helpful comments which significantly improved the presentation.

\printbibliography[title={References}]

\end{document}